\documentclass[11pt]{article}%

\usepackage{amsthm,amsmath,amsfonts}

\newtheorem{thm}{Theorem}[section]

\theoremstyle{definition}

\theoremstyle{remark}

\numberwithin{equation}{section}

\begin{document}

\title{A More General Maximal Bernstein-type Inequality}

\author{P\'eter Kevei
\thanks{Supported by the TAMOP--4.2.1/B--09/1/KONV--2010--0005 project.} \\
MTA-SZTE Analysis and Stochastics Research Group \\
Bolyai Institute, Aradi v\'ertan\'uk tere 1, 6720 Szeged, Hungary \\
e-mail: \texttt{kevei@math.u-szeged.hu} \smallskip \\ 
David M. Mason\thanks{Research partially supported by NSF
Grant DMS--0503908.}\\
University of Delaware \\
213 Townsend Hall, Newark, DE 19716, USA \\
e-mail: \texttt{davidm@udel.edu}}

\date{May 24, 2012}

\maketitle

\begin{abstract}
We extend a general
Bernstein-type maximal inequality of Kevei and Mason (2011) for sums of random variables.


\textit{AMS Subject Classificiation:} MSC 60E15; MSC 60F05; MSC 60G10.
\end{abstract}

\section{Introduction}

Let $X_{1},X_{2},\ldots$ be \ a sequence of random variables, and for any
choice of $1\leq k\leq l<\infty$ we denote the partial sum $S(k,l)=\sum
_{i=k}^{l}X_{i}$, and define $M(k,l)=\max\{|S(k,k)|,\ldots,|S(k,l)|\}$. It
turns out that under a variety of assumptions the partial sums $S(k,l)$ will
satisfy a generalized Bernstein-type inequality of the following form: for
suitable constants $A>0$, $a>0$, $b\geq0$ and $0<\gamma<2$ for all $m\geq0$,
$n\geq1$ and $t\geq0$,
\begin{equation}
\mathbf{P}\{|S(m+1,m+n)|>t\}\leq A\exp\left\{  -\frac{at^{2}}{n+bt^{\gamma}%
}\right\}  .\label{assumpgamma}%
\end{equation}

Kevei and Mason \cite{KM} provide numerous examples of sequences of random
variables $X_{1},X_{2},\ldots,$ that satisfy a Bernstein-type inequality of
the form (\ref{assumpgamma}). They show, somewhat unexpectedly, without any
additional assumptions, a modified version of it also holds for $M(1+m,n+m)$
for all $m\geq0$ and $n\geq1$. Here is their main result.

\begin{thm}
\label{T1} Assume that for constants $A>0$, $a>0$, $b\geq0$ and $\gamma
\in(0,2)$, inequality (\ref{assumpgamma}) holds for all $m\geq0,n\geq1$ and
$t\geq0$. Then for every $0<c<a$ there exists a $C>0$ depending only on $A,a$,
$b$ and $\gamma$ such that for all $n\geq1$, $m\geq0$ and $t\geq0$,
\begin{equation}
\mathbf{P}\{M(m+1,m+n)>t\}\leq C\exp\left\{  -\frac{ct^{2}}{n+bt^{\gamma}%
}\right\}  .\label{Cineq}%
\end{equation}
\end{thm}

There exists an interesting class of Bernstein-type inequalities that are not
of the form (\ref{assumpgamma}). Here are two motivating examples.\medskip

\noindent
\textbf{Example 1.}
Assume that $X_{1},X_{2},\dots,$ is a stationary
Markov chain satisfying the conditions of Theorem 6 of Adamczak \cite{Ad} and
let $f$ be any bounded measurable function such that $Ef\left(  X_{1}\right)
=0$. His theorem implies that for some constants $D>0$, $d_{1}>0 $ and
$d_{2}>0$ for all $t\geq0$ and $n\geq1$,%
\begin{equation}
P\left\{  |S_{n}(f)|\geq t\right\}  \leq D^{-1}\exp\left(  -\frac{Dt^{2}%
}{nd_{1}+td_{2}\log n}\right)  ,\label{exin1}%
\end{equation}
where $S_{n}(f)=\sum_{i=1}^{n}f\left(  X_{i}\right)$, and
$D/d_1$ is related to the limiting variance in the central limit theorem.
\smallskip

\noindent
\textbf{Example 2.}
Assume that $X_{1},X_{2},\dots,$ is a strong
mixing sequence with mixing coefficients $\alpha\left(  n\right)  $, $n\geq1$,
satisfying for some $d>0$, $\alpha\left(  n\right)  \leq\exp\left(
-2dn\right)  $. Also assume that $EX_{i}=0$ and for some $M>0$, $\left\vert
X_{i}\right\vert \leq M$, for all $i\geq1$. Theorem 2 of Merlev\`{e}de,
Peligrad and Rio \cite{MPR} implies that for some constant $D>0$ for all
$t\geq0$ and $n\geq1$,%
\begin{equation}
P\left\{  |S_{n}|\geq t\right\}  \leq D\exp\left(  -\frac{Dt^{2}}{nv^{2}%
+M^{2}+tM\left(  \log n\right)  ^{2}}\right)  ,\label{exin2}%
\end{equation}
where $S_{n}=\sum_{i=1}^{n}X_{i}$ and $v^{2}=\sup_{i>0}\left(  Var\left(
X_{i}\right)  +2\sum_{j>i}\left\vert cov\left(  X_{i},X_{j}\right)
\right\vert \right)  .$
\smallskip

The purpose of this note to establish the following extended version of
Theorem \ref{T1} that will show that a maximal version of inequalities
(\ref{exin1}) and (\ref{exin2}) also holds.

\begin{thm}
\label{T2} Assume that there exist constants $A>0$ and $a>0$ and a sequence of
non-decreasing non-negative functions $\left\{  g_{n}\right\}  _{n\geq1}$ on
$\left(  0,\infty\right)  $, such that for all $t>0$ and $n\geq1$,
$g_{n}\left(  t\right)  \leq g_{n+1}\left(  t\right)  $ and for all
$0<\rho<1$
\begin{equation}
\lim_{n\rightarrow\infty}\inf\left\{  \frac{t^{2}}{g_{n}(t)\log t}%
:g_{n}\left(  t\right)  >\rho n\right\}  =\infty,\label{function-assump}%
\end{equation}
where the infimum of the empty set is defined to be infinity, such that for all
$m\geq0$, $n\geq1$ and $t\geq0$,
\begin{equation}
P\{|S(m+1,m+n)|>t\}\leq A\exp\left\{  -\frac{at^{2}}{n+g_{n}(t)}\right\}
.\label{assumpg}%
\end{equation}
Then for every $0<c<a$ there exists a $C>0$ depending only on $A,a$ and
$\left\{  g_{n}\right\}  _{n\geq1}$ such that for all $n\geq1$, $m\geq0$ and
$t\geq0$,
\begin{equation}
P\{M(m+1,m+n)>t\}\leq C\exp\left\{  -\frac{ct^{2}}{n+g_{n}(t)}\right\}
.\label{max-ineq}%
\end{equation}
\end{thm}

Note that condition (\ref{function-assump}) trivially holds when the functions
$g_n$ are bounded, since the corresponding sets are empty sets.
However, in the interesting cases $g_n$'s are not bounded, and in this case
the condition basically says that $g_n(t)$ increases slower than $t^2$.

Essentially the same proof shows that the statement of Theorem \ref{T2} remains
true if in the numerator of (\ref{assumpg}) and (\ref{max-ineq})
the function $t^2$ is replaced by a regularly
varying function at infinity $f(t)$ with a positive index. In this case the $t^2$
in condition (\ref{function-assump}) must be replaced by $f(t)$.
Since we do not know any application of a result of this type, we only mention this 
generalization.
\begin{proof}
Choose any $0<c<a.$ We prove our theorem by induction on $n$. Notice
that by the assumption, for any integer $n_{0}\geq1$ we may choose $C>An_{0}$
to make the statement true for all $1\leq n\leq n_{0}$. This remark will be
important, because at some steps of the proof we assume that $n$ is large
enough. Also since the constants $A$ and $a$ in (\ref{assumpg}) are
independent of $m$, we can without loss of generality assume $m=0$.

Assume the statement holds up to some $n\geq2$. (The constant $C$ will be
determined in the course of the proof.)

\noindent\textbf{Case 1.} Fix a $t>0$ and assume that 
\begin{equation}
g_{n+1}(t)\leq\alpha\,n,\label{alpha}%
\end{equation}
for some $0<\alpha<1$ be specified later. (In any case, we assume that $\alpha
n\geq1$.) Using an idea of \cite{MSS}, we may write for arbitrary $1\leq k<n$,
$0 < q < 1$ and $p+q=1$ the inequality%
\begin{align*}
P\{M(1,n+1)>t\}\leq &  P\{M(1,k)>t\}+P\{|S(1,k+1)|>pt\}\\
&  +P\{M(k+2,n+1)>qt\}.
\end{align*}

Let%
\[
u=\frac{n+g_{n+1}(qt)-q^{2}g_{n+1}(t)}{1+q^{2}}.
\]
Note that $u \leq n-1$ if $0 < \alpha < 1$ is chosen small enough depending
on $q$, for $n$ large enough.
Notice that
\begin{equation}
\frac{t^{2}}{u+g_{n+1}(t)}=\frac{q^{2}t^{2}}{n-u+g_{n+1}(qt)}.\label{eq}%
\end{equation}
Set
\begin{equation}
k=\left\lceil u\right\rceil .\label{k}%
\end{equation}
Using the induction hypothesis and (\ref{assumpg}), keeping in mind that
$1 \leq k \leq n-1$, we obtain%
\begin{equation}
\begin{split}
P\{M(1,n+1)>t\}
\leq & C\exp\left\{  -\frac{ct^{2}}{k+g_{k}(t)}\right\}
+A\exp\left\{  -\frac{ap^{2}t^{2}}{k+1+g_{k+1}(pt)}\right\}  \\
&  +C\exp\left\{  -\frac{cq^{2}t^{2}}{n-k+g_{n-k}(qt)}\right\} \\
\leq & C\exp\left\{  -\frac{ct^{2}}{k+g_{n+1}(t)}\right\}  +A\exp\left\{
-\frac{ap^{2}t^{2}}{k+1+g_{n+1}(pt)}\right\} \\
& +C\exp\left\{  -\frac
{cq^{2}t^{2}}{n-k+g_{n+1}(qt)}\right\}. \label{mainineqgamma}%
\end{split}
\end{equation}
Notice that we chose $k$ to make the first and third terms in
(\ref{mainineqgamma}) almost equal, and since by (\ref{k})
\[
\frac{t^{2}}{k+g_{n+1}(t)}\leq\frac{q^{2}t^{2}}{n-k+g_{n+1}(qt)}%
\]
the first term is greater than or equal to the third.

First we handle the second term in formula (\ref{mainineqgamma}), showing that whenever
$g_{n+1}(t)\leq\alpha n$,
\[
\exp\left\{  -\frac{ap^{2}t^{2}}{k+1+g_{n+1}(pt)}\right\}  \leq\exp\left\{
-\frac{ct^{2}}{n+1+g_{n+1}(t)}\right\}  .
\]
For this we need to verify that for $g_{n+1}(t)\leq\alpha n$,
\begin{equation}
\frac{ap^{2}}{k+1+g_{n+1}(pt)}>\frac{c}{n+1+g_{n+1}(t)},\label{eq1}%
\end{equation}
which is equivalent to
\[
ap^{2}(n+1+g_{n+1}(t))>c(k+1+g_{n+1}(pt)).
\]
Using that
\[
k=\lceil u\rceil\leq u+1=1+\frac{1}{1+q^{2}}\left[  n+g_{n+1}(qt)-q^{2}%
g_{n+1}(t)\right]  ,
\]
it is enough to show%
\[
n\left(  ap^{2}-\frac{c}{1+q^{2}}\right)  +ap^{2}-2c
\]
\[
+\left[  g_{n+1}(t)ap^{2}-g_{n+1}(pt)c-\frac{c}{1+q^{2}}\left(  g_{n+1}%
(qt)-q^{2}g_{n+1}(t)\right)  \right]  >0.
\]
Note that if the coefficient of $n$ is positive, then we can choose $\alpha$
in (\ref{alpha}) small enough to make the above inequality hold. So in order
to guarantee (\ref{eq1}) (at least for large $n$) we only have to choose the
parameter $p$ so that $ap^{2}-c>0$, which implies that
\begin{equation}
ap^{2}-\frac{c}{1+q^{2}}>0\label{assump-c-1}%
\end{equation}
holds, and then select $\alpha$ small enough, keeping mind that we assume
$\alpha n \geq 1$ and $k \leq n-1$.

Next we treat the first and third terms in (\ref{mainineqgamma}). Because of
the remark above, it is enough to handle the first term. Let us examine the
ratio of $C\exp\{-ct^{2}/(k+g_{n+1}(t))\}$ and $C\exp\{-ct^{2}/(n+1+g_{n+1}%
(t))\}$. Notice again that since $u+1\geq k$, the monotonicity of $g_{n+1}(t)$
and $g_{n+1}(t)\leq\alpha n$ implies
\begin{align*}
n+1-k  &  \geq n-u=n-\frac{n+g_{n+1}(qt)-q^{2}g_{n+1}(t)}{1+q^{2}}\\
&  \geq\frac{q^{2}n-(1-q^{2})g_{n+1}(t)}{1+q^{2}}\\
&  \geq n\frac{q^{2}-\alpha(1-q^{2})}{1+q^{2}}\\
&  =:c_{1}n.
\end{align*}
At this point we need that $0<c_{1}<1$. Thus we choose $\alpha$ small enough
so that
\begin{equation}
q^{2}-\alpha(1-q^{2})>0.\label{assump-c-2}%
\end{equation}
Also we get using $g_{n+1}(t)\leq\alpha n$ the bound
\[
(n+1+g_{n+1}(t))(k+g_{n+1}(t))\leq2n^{2}(1+\alpha)^{2}=:c_{2}n^{2},
\]
which holds if $n$ large enough. Therefore, we obtain for the ratio
\[
\exp\left\{  -ct^{2}\left(  \frac{1}{k+g_{n+1}(t)}-\frac{1}{n+1+g_{n+1}%
(t)}\right)  \right\}  \leq\exp\left\{  -\frac{cc_{1}t^{2}}{c_{2}n}\right\}
\leq\mathrm{e}^{-1},
\]
whenever $cc_{1}t^{2}/(c_{2}n)\geq1$, that is $t\geq\sqrt{c_{2}n/(cc_{1})}$.
Substituting back into (\ref{mainineqgamma}), for $t\geq\sqrt{c_{2}n/(cc_{1}%
)}$ and $g_{n+1}(t)\leq\alpha n$ we obtain
\[
P\{M(1,n+1)>t\}
\]%
\[
\leq\left(  \frac{2}{\mathrm{e}}C+A\right)  \exp\{-ct^{2}/(n+1+g_{n+1}%
(t))\}\leq C\exp\{-ct^{2}/(n+1+g_{n+1}(t))\},
\]
where the last inequality holds for $C>A\mathrm{e}/(\mathrm{e}-2)$.

Next assume that $t<\sqrt{c_{2}n/(cc_{1})}$. In this case choosing $C$ large
enough we can make the bound $>1$, namely
\[
C\exp\left\{  -\frac{ct^{2}}{n+1+g_{n+1}(t)}\right\}  \geq C\exp\left\{
-\frac{cc_{2}n}{cc_{1}n}\right\}  =C\mathrm{e}^{-c_{2}/c_{1}}\geq1,
\]
if $C>\mathrm{e}^{c_{2}/c_{1}}$.\medskip

\noindent\textbf{Case 2.} Now we must handle the case $g_{n+1}(t)>\alpha n$.
Here we apply the inequality
\[
P\{M(1,n+1)>t\}\leq P\{M(1,n)>t\}+P\{|S(1,n+1)|>t\}.
\]
Using assumption (\ref{assumpg}) and the induction hypothesis, we have
\begin{equation*}
\begin{split}
P\{M(1,n+1)>t\}
& \leq C\exp\left\{  -\frac{ct^{2}}{n+g_{n}(t)}\right\}
+A\exp\left\{  -\frac{at^{2}}{n+1+g_{n+1}(t)}\right\} \\
& \leq C\exp\left\{  -\frac{ct^{2}}{n+g_{n+1}(t)}\right\}  +A\exp\left\{
-\frac{at^{2}}{n+1+g_{n+1}(t)}\right\}  .
\end{split}
\end{equation*}
We will show that the right side $\leq C\exp\{-ct^{2}/(n+1+g_{n+1}(t))\}$. For
this it is enough to prove
\begin{equation}
\begin{split}
& \exp\left\{  -ct^{2}\left(  \frac{1}{n+g_{n+1}(t)}-\frac{1}{n+1+g_{n+1}%
(t)}\right)  \right\} \\
& +\frac{A}{C}\exp\left\{  -\frac{t^{2}(a-c)}%
{n+1+g_{n+1}(t)}\right\}  \leq1.\label{less}%
\end{split}
\end{equation}
Using the bound following from $g_{n+1}(t)>\alpha n$ and recalling that
$\alpha n\geq1$ and $0 < \alpha < 1$, we get
\[
\frac{t^{2}}{(n+g_{n+1}(t))(n+1+g_{n+1}(t))}\geq\frac{ \alpha^2 t^{2}}
{(1 + \alpha)(1 + 2 \alpha)g_{n+1}(t)^{2}}=:c_{3}\frac{t^{2}}{g_{n+1}(t)^{2}},
\]
and
\[
\frac{t^{2}(a-c)}{n+1+g_{n+1}(t)}\geq\frac{t^{2}}{g_{n+1}(t)}\frac
{\alpha (a-c)}{1 + 2 \alpha}=:\frac{t^{2}}{g_{n+1}(t)}c_{4}.
\]

Choose $\delta> 0$ so small such that $0 < x \leq\delta$ implies
$\mathrm{e}^{-c c_{3} x^{2}} \leq1 - \frac{c c_{3}}{2} x^{2}$.

\noindent
For $t/g_{n+1}(t)\geq\delta$ the left-hand side of (\ref{less}) is less then
\[
\mathrm{e}^{-cc_{3}\delta^{2}}+\frac{A}{C},
\]
which is less than 1, for $C$ large enough.

For $t/g_{n+1}(t)\leq\delta$ by the choice of $\delta$ the left-hand side of
(\ref{less}) is less then
\[
1-\frac{cc_{3}}{2}\frac{t^{2}}{g_{n+1}(t)^{2}}+\frac{A}{C}\exp\left\{
-\frac{t^{2}}{g_{n+1}(t)}c_{4}\right\}  ,
\]
which is less than 1 if
\[
\frac{cc_{3}}{2}\frac{t^{2}}{g_{n+1}(t)^{2}}>\frac{A}{C}\exp\left\{
-\frac{t^{2}}{g_{n+1}(t)}c_{4}\right\}  .
\]
By (\ref{function-assump}), for any $0<\eta<1$ and all large enough $n$,
$g_{n+1}(t)1\left\{  g_{n+1}\left(  t\right)  >\alpha n\right\}  \leq\eta
t^{2}$, so that for all large $n,$ whenever $g_{n+1}\left(  t\right)  >\alpha
n,$ we have
\[
\frac{t^{2}}{g_{n+1}(t)^{2}}\geq t^{-2},
\]
and again by (\ref{function-assump}) for all large $n,$ whenever
$g_{n+1}\left(  t\right)  >\alpha n,$ $t^{2}/g_{n+1}(t)\geq\left(  3/c_{4}\,\right)
\log t$. Therefore for all large $n,$ whenever $g_{n+1}\left(  t\right)
\alpha n,$
\[
\exp\left\{  -\frac{t^{2}}{g_{n+1}(t)}c_{4}\right\}  \leq t^{-3},
\]
which is smaller than $t^{-2}\frac{Ccc_{3}}{2A}$, for $t$ large enough, i.e.
for $n$ large enough. The proof is complete.
\end{proof}

By choosing $g_{n}\left(  t\right)  =bt^{\gamma}$ for all $n\geq1$ we see that
Theorem \ref{T2} gives Theorem \ref{T1} as a special case. Also note that
Theorem \ref{T2} remains valid for sums of Banach space valued random
variables with absolute value $\left\vert \cdot\right\vert $ replaced by norm
$\left\vert |\cdot|\right\vert $. Theorem \ref{T2} permits us to derive the
following maximal versions of inequalities (\ref{exin1}) and (\ref{exin2}).
\medskip

\noindent\textbf{Application 1}. In Example 1 one readily checks that the
assumptions of Theorem \ref{T2} are satisfied with $A=D^{-1}$ and $a=D/d_{1}$
\[
g_{n}\left(  t\right)  =\left(  \frac{td_{2}}{d_{1}}\right)  \log n\text{.}%
\]
We get the maximal version of inequality (\ref{exin1}) holding for any $0<c<1
$ and all $n\geq1$ and $t>0$
\begin{equation}
P\left\{  \big|\max_{1\leq m\leq n}S_{n}(f) \big| \geq t\right\}  \leq C\exp\left(
-\frac{cDt^{2}}{nd_{1}+td_{2}\log n}\right)  ,\label{MP1}%
\end{equation}
for some constant $C\geq D^{-1}$ depending on $c$, $D^{-1}$, $D/d_{1}$ and
$\left\{  g_{n}\right\}  _{n\geq1}$.\smallskip

\noindent\textbf{Application 2.} In Example 2 one can verify that the
assumptions of the Theorem \ref{T2} hold with $A=D$ and $a=D/v^{2}$ and
\[
g_{n}\left(  t\right)  =\frac{M^{2}}{v^{2}}+\left(  \frac{tM}{v^{2}}\right)
\left(  \log n\right)  ^{2}\text{,}%
\]
which leads to the maximal version of inequality (\ref{exin2}) valid for any
$0<c<1$ and all $n\geq1$ and $t>0$
\begin{equation}
P\left\{  \max_{1\leq m\leq n}|S_{m}|\geq t\right\}  \leq C\exp\left(
-\frac{cDt^{2}}{nv^{2}+M^{2}+tM\left(  \log n\right)  ^{2}}\right) \label{MP}%
\end{equation}
for some constant $C\geq D$ depending on $c$, $D/v^{2}$ and $\left\{
g_{n}\right\}  _{n\geq1}$. See Corollary 24 of Merlev\`ede and Peligrad
\cite{MP} for a closely related inequality that holds for all $n\geq2$ and
$t>K\log n$ for some $K>0.\medskip$

\noindent\textbf{Remark} There is a small oversight in the published version
of the Kevei and Mason paper. Here are the corrections that fix it.\smallskip

\noindent1. Page 1057, line -9: Replace \textquotedblleft$1\leq k\leq
n$\textquotedblright\ by \textquotedblleft$1\leq k<n$\textquotedblright.

\noindent2. Page 1057, line -7: Replace this line with

$\leq\mathbf{P}\left\{  M\left(  1,k\right)  >t\right\}  +\mathbf{P}\left\{
S\left(  1,k+1\right)  >pt\right\}  +\mathbf{P}\left\{  M\left(
k+2,n+1\right)  >qt\right\}  .$

\noindent3. Page 1058: Replace \textquotedblleft$k+bp^{\gamma}t^{\gamma}%
$\textquotedblright\ by \textquotedblleft$k+1+bp^{\gamma}t^{\gamma}%
$\textquotedblright\ in equations (2.4) and (2.5), as well as in line -13.

\noindent4. Page 1058: Replace \textquotedblleft$ap^{2}-c$\textquotedblright%
\ by \textquotedblleft$ap^{2}-2c$\textquotedblright\ in line -9.

\subsection*{Acknowledgment}

We thank a referee for a careful reading of the manuscript and a number of useful
comments.

% ------------------------------------------------------------------------
\end{document}